\documentclass[journal]{IEEEtran}
\usepackage{amsmath}
\usepackage{amssymb}
\usepackage{graphicx}
\usepackage{color}
\usepackage{subcaption}

\newtheorem{rem}{Remark}
\newtheorem{defn}{Definition}
\newtheorem{prop}{Proposition}

\newtheorem{lem}{Lemma}
\newtheorem{alg}{Algorithm}

\newcommand{\NA}{N_A}
\newcommand{\ND}{N_D}

\usepackage{my_macros}

\title{\LARGE Multiplayer Reach-Avoid Games via Pairwise Outcomes}

\author{Mo Chen, Zhengyuan Zhou and Claire J. Tomlin%
\thanks{This work has been supported in part by NSF under CPS:ActionWebs (CNS-931843), by ONR under the HUNT (N0014-08-0696) and SMARTS (N00014-09-1-1051) MURIs and by grant N00014-12-1-0609, by AFOSR under the CHASE MURI (FA9550-10-1-0567). Recommended by Senior Editor Edwin K. P. Chong.}
\thanks{M.~Chen, and C.~J.~Tomlin are with the Department of Electrical Engineering and Computer Sciences,
        University of California, Berkeley, CA 94720, USA
        {\tt\small \{mochen72,tomlin\}@eecs.berkeley.edu}}
\thanks{Z.~Zhou is with the Department of Electrical Engineering,
        Stanford University, Stanford, CA 94305, USA
        {\tt\small zyzhou@stanford.edu}}   
\thanks{We thank Haomiao Huang for sharing MATLAB code for 4D HJI calculations.}
}

\begin{document}

\maketitle
\thispagestyle{empty}
\pagestyle{empty}

\begin{abstract}
A multiplayer reach-avoid game is a differential game between an attacking team with $\NA$ attackers and a defending team with $\ND$ defenders playing on a compact domain with obstacles. The attacking team aims to send $\m$ of the $\NA$ attackers to some target location, while the defending team aims to prevent that by capturing attackers or indefinitely delaying attackers from reaching the target. Although the analysis of this game plays an important role in many applications, the optimal solution to this game is computationally intractable when $\NA>1$ or $\ND>1$. In this paper, we present two approaches for the $\NA=\ND=1$ case to determine pairwise outcomes, and a graph theoretic maximum matching approach to merge these pairwise outcomes for an $\NA,\ND>1$ solution that provides guarantees on the performance of the defending team. We will show that the four-dimensional Hamilton-Jacobi-Isaacs approach allows for real-time updates to the maximum matching, and that the two-dimensional ``path defense" approach is considerably more scalable with the number of players while maintaining defender performance guarantees.
\end{abstract}


\section{Introduction}

Multiplayer reach-avoid games are differential games between two adversarial teams of cooperative players playing on a compact domain with obstacles. The ``attacking team'' aims to send as many team members, called ``attackers'', to some target set as quickly as possible. The ``defending team'' seeks to delay or prevent the attacking team from doing so by attempting to capture the attackers. Such differential games have been studied extensively \cite{HThesis, Huang2011} and are also powerful theoretical tools for analyzing realistic situations in robotics, aircraft control, security, and other domains \cite{OFTHEAIRFORCEWASHINGTON:2009p37, Erzberger:2006p44, kiva2009}.


The multiplayer reach-avoid game is difficult to analyze because the two teams have conflicting and asymmetric goals while complex cooperation within each team may exist. In addition, optimal solutions are impossible to compute using traditional dynamic programming approaches due to the intrinsic high dimensionality of the joint state space. Previously, in \cite{Earl:2007p101}, where a team of defenders assumes that the attackers move towards their target in straight lines, a mixed-integer linear programming approach was used. \cite{Chasparis:2005p102} assumes that the attackers use a linear feedback control law, and a mixed integer linear program was then relaxed into linear program. In complex pursuit-evasion games where players may change roles over time, a nonlinear model-predictive control \cite{Sprinkle:2004p100} approach has been investigated. Approximate dynamic programming \cite{McGrew:2008p103} has also been used to analyze reach-avoid games.

Although the above techniques provide some useful insight, they only work well when strong assumptions are made or when accurate models of the opposing team can be obtained. To solve general reach-avoid games, the Hamilton-Jacobi-Isaacs (HJI) approach \cite{b:isaacs-1967} is ideal when the game is low-dimensional. The approach involves solving an HJI partial differential equation (PDE) in the joint state space of the players to compute a reach-avoid set, which partitions the players' joint state space into a winning region for the defending team and one for the attacking team. The optimal strategies can then be extracted from the gradient of the solution. This approach is particularly useful because of the numerical tools \cite{j:mitchell-TAC-2005, Sethian1996, b:osher-fedkiw-2002} available, and has been able to solve several practical problems \cite{Huang2011, j:mitchell-TAC-2005, DSST08}. The HJI approach can be applied to a large variety of player dynamics and does not explicitly assume any control strategy or prediction models for the players. However, the approach cannot be directly applied to our multiplayer reach-avoid game because its complexity scales exponentially with the number of players, making the approach only tractable for the two-player game. Thus, complexity-optimality trade-offs must be made.


For the two-player reach-avoid game, we first present the two-player HJI solution \cite{Huang2011}, which computes a four-dimensional (4D) reach-avoid set that determines which player wins the game assuming both players use the closed-loop optimal control strategy. Next, we propose the ``path defense" approximation to the HJI solution, in which the defenders utilize a ``semi-open-loop" control strategy. Here, we approximate two-dimensional (2D) slices of the reach-avoid sets by solving 2D Eikonal equations, and provide guarantees for the defending team's performance.

For the multiplayer reach-avoid game, we propose to merge the $\NA\ND$ pairwise outcomes using the graph theoretic maximum matching, which can be efficiently computed by known algorithms \cite{Schrjiver2004, Karpinski1998}. The maximum matching process incorporates cooperation among defenders without introducing significant additional computation cost. When players on each team have identical dynamics, only a \textit{single} HJI PDE needs to be solved to characterize \textit{all} pairwise outcomes. Furthermore, when applying maximum matching to the two-player path defense solution, the computational complexity scales linearly with the number of attackers, as opposed to quadratically with the total number of players in the HJI approach. 


\section{The Reach-Avoid Problem}
\subsection{The Multiplayer Reach-Avoid Game}
\label{sec:formulation}
Consider $\NA+\ND$ players partitioned into the set of $\NA$ attackers, $\pas = \{\pam{1}, \pam{2}, \ldots, \pam{\NA}\}$ and the set of $\ND$ defenders, $\pbs = \{\pbm{1}, \ldots, \pbm{\ND}\}$, whose states are confined in a bounded, open domain $\amb \subset \R^2$. The domain $\amb$ is partitioned into $\amb$ = $\free \cup \obs$, where $\free$ is a compact set representing the free space in which the players can move, while $\obs = \amb \setminus \free$ corresponds to obstacles in the domain. 

Let $\xam{i}, \xbm{j} \in \R^2$ denote the state of players $\pam{i}$ and $\pbm{j}$, respectively. Then given initial conditions $\xanm{i}\in \free,i=1,2,\ldots,\NA,\xbnm{i}\in \free,i=1,2,\ldots,\ND$, we assume the dynamics of the players to be defined by the following decoupled system for $t \geq 0$:

\bq\label{eq:dynamics}
\begin{aligned}
\dotxam{i}(t) &= \velai{i}\cam{i}(t), & \xam{i}(0) = \xanm{i}, i=1,2,\ldots,\NA \\
\dotxbm{i}(t) &= \velbi{i}\cbm{i}(t), & \xbm{i}(0) = \xbnm{i}, i=1,2,\ldots,\ND
\end{aligned}
\eq
where $\velai{i}, \velbi{i}$ denote maximum speeds for $\pam{i}$ and $\pbm{i}$ respectively, and $\cam{i},\cbm{i}$ denote controls of $\pam{i}$ and $\pbm{i}$ respectively. We assume that $\cam{i},\cbm{i}$ are drawn from the set $\A = \{\sigma \colon [0,\infty)\rightarrow \unitball \mid \sigma \text{ is measurable}\}$, where $\unitball$ denotes the closed unit disk in $\R^2$. We also constrain the players to remain within $\free$ for all time. Denote the joint state of all players by $\xj = (\xja, \xjb)$ where $\xja =(\xam{1},\ldots\xam{\NA})$ is the attacker joint state $\pas$, and $\xjb = (\xbm{1},\ldots,\xbm{\ND})$ is the defender joint state $\pbs$. 

\begin{figure}
\centering
\includegraphics[width=0.35\textwidth]{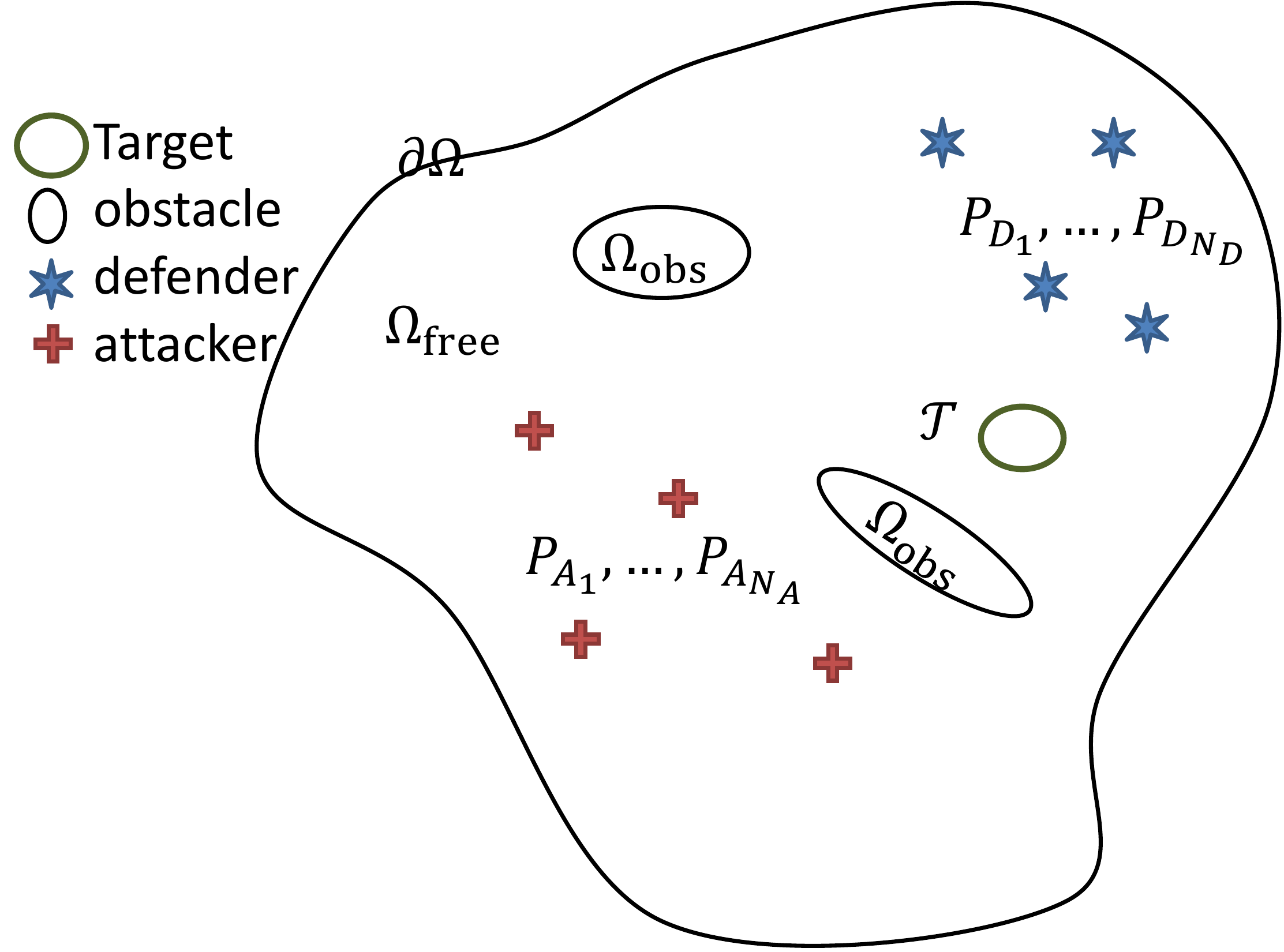}
\caption{The components of a multiplayer reach-avoid game.}
\label{fig:mp_form}
\end{figure}

The attacking team wins whenever $\m$ of the $\NA$ attackers reach some target set without being captured by the defenders; $\m$ is pre-specified with $0<M\le \NA$. The target set is denoted $\target\subset\free$ and is compact. The defending team wins if it can prevent the attacking team from winning by capturing or indefinitely delaying $\NA-\m+1$ attackers from reaching $\target$. An illustration of the game setup is shown in Fig. \ref{fig:mp_form}.

Let $\avoid_{ij} = \left\{\xj\in\amb^{\NA+\ND} \mid \|\xam{i}-\xbm{j}\|_2\le\Rc \right\}$ denote the capture set. $\pam{i}$ is captured by $\pbm{j}$ if $\pam{i}$'s position is within a distance $\Rc$ of $\pbm{j}$'s position. 

In this paper, we address the following problems:
\begin{enumerate}
\item Given $\xjn$, $\target$, and some fixed integer $\m, 0<\m\le\NA$, can the attacking team win?
\item More generally, given $\xjn$ and $\target$, how many attackers can the defending team prevent from reaching the target?
\end{enumerate}

\subsection{The Two-Player Reach-Avoid Game}
\label{sec:2p_ra}
We will answer the above questions about the $\NA$ vs. $\ND$ reach-avoid game by using the solution to the two-player $1$ vs. $1$ game as a building block. In the two-player game, we denote the attacker $\pa$, the defender $\pb$, their states $\xa,\xb$, and their initial conditions $\xan,\xbn$. Their dynamics are
\bq
\begin{aligned}
\dotxa(t) &= \vela\ca(t), & \xa(0) = \xan,\\
\dotxb(t) &= \velb\cb(t), & \xb(0) = \xbn
\end{aligned}
\eq

The players' joint state becomes $\xj=(\xa,\xb)$, and their joint initial condition becomes $\xjn=(\xan,\xbn)$. The capture set becomes simply $\avoid = \left\{(\xa,\xb)\in\amb^2 \mid \|\xa-\xb\|_2\leq \Rc\right\}$. 

$\pa$ wins if it reaches the target $\target$ without being captured by $\pb$. $\pb$ wins if it can prevent $\pa$ from winning by capturing $\pa$ or indefinitely delaying $\pa$ from reaching $\target$. For the two-player reach-avoid game, we seek to answer the following:
\begin{enumerate}
\item Given $\xjn$ and $\target$, is the defender guaranteed to win? \label{p:tp1}
\item More generally, given $\xa$ and $\target$, what is the set of initial positions from which the defender is guaranteed to win? \label{p:tp2}
\end{enumerate}

\section{The HJI Solution of the 1 vs. 1 Game} \label{sec:solution_hji}
The HJI approach for solving differential games is outlined in \cite{Huang2011,j:mitchell-TAC-2005, LSToolbox}. The optimal joint closed-loop control strategies for the attacker and the defender in a two-player reach-avoid game can be obtained by solving a 4D HJI PDE. This solution allows us to determine whether the defender will win against the attacker in a 1 vs. 1 setting. 

In the two-player game, the attacker aims to reach $\target$ while avoiding $\avoid$. Both players also avoid $\obs$. In particular, the defender wins if the attacker is in $\obs$, and vice versa. Therefore, we define the terminal set and avoid set to be

\begin{equation} \label{eq:4DHJI_sets}
\begin{aligned}
R &= \left\{\xj\in\amb^2 \mid \xa\in\target \right\} \cup \left\{\xj\in\amb^2\mid \xb\in\obs \right\} \\
A &= \avoid \cup \left\{\xj\in\amb^2\mid \xa\in\obs \right\}
\end{aligned} 
\end{equation}

Given \eqref{eq:4DHJI_sets}, we can define the corresponding implicit surface functions $\valsR_R,\valsR_A$ required for solving the HJI PDE. Since $\amb\subset\R^2$, the result is $\mathcal{RA}_\infty(R,A)\in\R^4$, a 4D reach-avoid set. If $\xjn\in\mathcal{RA}_\infty(R,A)$, then the attacker is guaranteed to win the game by using the optimal control \textit{even if} the defender is also using the optimal control; if $\xjn\notin\mathcal{RA}_\infty(R,A)$, then the defender is guaranteed to win the game by using the optimal control \textit{even if} the attacker is also using the optimal control.

\section{The Path Defense Solution to the 1 vs. 1 Game}
\label{sec:path_defense}
We approximate 2D slices of the 4D reach-avoid set (or simply ``2D slices") in the path defense approach. Each slice will be taken at an attacker position. Here, we will assume that the defender is not slower than the attacker: $\vela \leq \velb$. 

\subsection{The Path Defense Game}
\label{subsec:pd_game}
\begin{figure}
\centering
\includegraphics[width=0.35\textwidth]{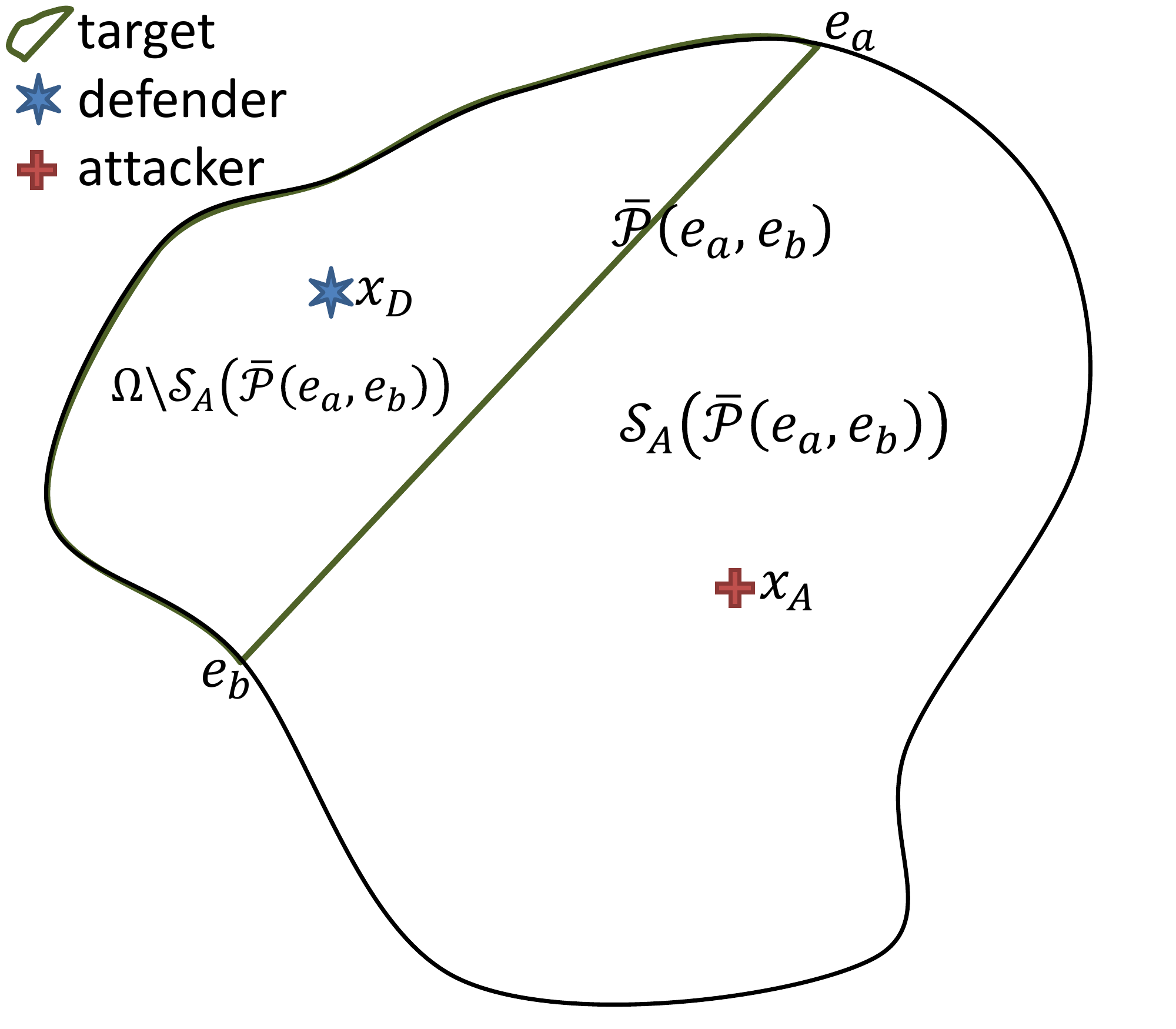}
\caption{The components of a path defense game.}
\label{fig:pd_form}
\end{figure}

The \textit{Path Defense Game} is a two-player reach-avoid game in which the boundary of the target set is the shortest path between two points on $\boundary$, and the target set is on one side of that shortest path (Fig. \ref{fig:pd_form}). We denote the target set as $\target=\sac$ for two given points on the boundary $\apa,\apb$. $\sa$ and $Anal.d$ are defined below. 

\begin{defn} 
\textbf{Path of defense}. A path of defense, $Anal.d$, is the shortest path between two boundary points $\apa,\apb\in\boundary$. $\apa$ and $\apb$ are called the \textbf{anchor points} of path $\pathd$. 
\end{defn}

Denote the shortest path between \textit{any} two points $\x,\y\in\free$ to be $\spath(\x,\y)$, with length $\dist(\x,\y)$, and requiring the attacker and defender durations of $\ta(\x, \y),\tb(\x,\y)$ to traverse, respectively. We will also use $\dist(\cdot,\cdot)$ with one or both arguments being sets in $\amb$ to denote the shortest distance between the arguments.

\begin{defn} 
\textbf{Attacker's side of the path}. A path of defense $\pathd$ partitions the domain $\amb$ into two regions. Define $\sa$ to be the region that contains the attacker, not including points on the path $\pathd$. The attacker seeks to reach the target set $\target=\sac$.
\end{defn}

\subsection{Solving The Path Defense Game}
A path defense game can be directly solved by computing a 4D reach-avoid set. Since the direct solution is time- and memory- intensive, we propose an efficient approximation of 2D slices that is conservative towards the defender.

\begin{defn} 
\textbf{Defendable path}. Given $\xjn=(\xan,\xbn)$, a path $\pathd$ is defendable if regardless of the attacker's actions, the defender has a control function $\cb(\cdot)$ to prevent the attacker from reaching $\pathd$ without being captured.
\end{defn}

\begin{defn} 
\textbf{Strongly defendable path}. $\pathd$ is \textit{strongly} defendable if regardless of the attacker's actions, the defender has a control function $\cb(\cdot)$ to reach $\pathd$ after finite time and prevent the attacker from reaching $\pathd$.
\end{defn}

\begin{figure}
\centering
\includegraphics[width=0.45\textwidth]{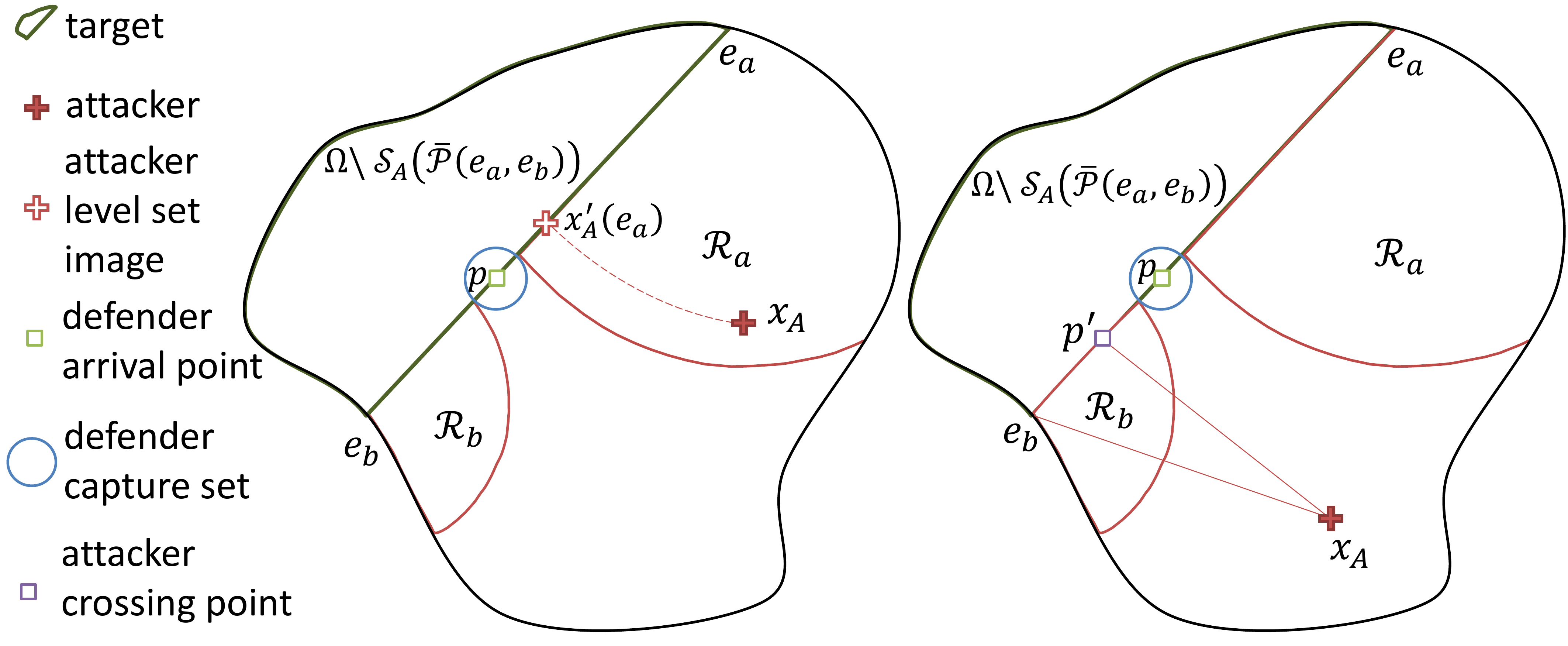}
\caption{Left: If the attacker is in $\rpa\cup\rpb$ and moves towards $e_a$, he  will be able to reach $\apa$ without being captured. Right: If the attacker is not in $\rpa\cup\rpb$, there is no point on the path $\pprime\in\pathd$ that the he can reach without being captured.}
\label{fig:lemma1}
\end{figure}

Checking whether a path $\pathd$ is defendable involves a 4D reach-avoid set calculation, so instead we check whether a path $\pathd$ is strongly defendable. The following definitions lead to our first Lemma which describes how to determine strong defendability using 2D distance calculations; the definitions and Lemma are illustrated in Fig. \ref{fig:lemma1}.

\begin{defn} 
\textbf{Attacker level set image}. Given attacker position $\xa(t)$, define the attacker level set image with respect to anchor point $\apa$ to be $\xai(t;\apa) = \{\x\in\pathd: \ta(\x,\apa) = \ta(\xa(t),\apa)\}$. $\xai$ is the unique point on $\pathd$ such that $\ta(\xai, \apa)=\ta(\xa,\apa)$. Define $\xai(t;\apb)$ similarly by replacing $\apa$ with $\apb$. For convenience, we sometimes omit the time argument and write $\xai(\apa)$.
\end{defn}

\begin{prop}
\label{rem:image_of_a}
$\dist(\xai(\apb),\apa) \le \dist(\xai(\apa),\apa)$. 
\end{prop}

\begin{IEEEproof}
First note that
\begin{equation*}
\begin{aligned}
\dist(\apa,\apb) &\le \dist(\xa,\apa) + \dist(\xa,\apb)\\
& = \dist(\xai(\apa),\apa) + \dist(\xai(\apb),\apb)
\end{aligned}
\end{equation*}

Then, since the left hand side is given by $\dist(\apa,\apb) = \dist(\apa,\xai(\apb))+\dist(\xai(\apb),\apb)$, the result follows.
\end{IEEEproof}

\begin{defn} 
\textbf{Capture set}: Define the capture set to be $\Dc(\y,t)=\{\x \mid \|\x-\y(t)\|_2\leq \Rc\}$. We will drop the second argument of $\Dc$ when $\y$ does not depend on time.
\end{defn}

\begin{rem}
Given $\pathd$, suppose the attacker level set image is within defender's capture set at some time $s$: 
\begin{equation*}
\begin{aligned}
&\xai(s;\apa) \in \Dc(\xb, s) \quad \text{(or } & \xai(s;\apb) \in \Dc(\xb, s) \text{)}
\end{aligned}
\end{equation*}
Then, there exists a control for the defender to keep the attacker level set image within the capture radius of the defender thereafter: 

\begin{equation*}
\begin{aligned}
&\xai(t;\apa) \in \Dc(\xb, t) ~ \forall t\geq s \\
\text{(or } & \xai(t;\apb) \in \Dc(\xb, t) ~ \forall t\geq s\text{)}
\end{aligned}
\end{equation*}

This is because the attacker level set image can move at most as fast as the attacker, who is not faster than the defender.
\end{rem}

\begin{defn} 
\label{def:d_win_region}
\textbf{Regions induced by point $\ppath$ on path}. Given a point $\ppath\in\pathd$, define a region $\rpa\left(\ppath\right)$ associated the point $\ppath$ and anchor point $\apa$ as follows:
\bq
\rpa\left(\ppath\right) = \left\{\x: \dist(\x,\apa) \leq \dist(\Dc(\ppath),\apa) \right\}
\eq

Define $\rpb(\ppath)$ similarly by replacing $\apa$ with $\apb$.
\end{defn}

\begin{lem}
\label{lem:d_winning_region}
Suppose $\xbn = \ppath \in \pathd$, and $\vela=\velb$. Then, $\pathd$ is strongly defendable if and only if $\xan$ is outside the region induced by $\ppath$: $\xan\in\amb\backslash\left(\rpa \cup \rpb\right)$.
\end{lem}

\begin{IEEEproof} 
See Fig. \ref{fig:lemma1}. Assume $\xan \notin \target = \sac$, otherwise the attacker would start inside the target set. 

First, we show that if $\xan\in \rpa \cup \rpb$, then the attacker can reach $\apa$ or $\apb$ and hence $\sac$ without being captured. Without loss of generality (WLOG), suppose $\xan\in\rpa$. To capture the attacker, the defender's capture set must contain $\xai(\apa)$ or $\xai(\apb)$ at some time $t$. By Definition \ref{def:d_win_region}, we have $\dist(\xan,\apa) < \dist(\Dc(\ppath),\apa)$, so $\ta(\xai(\apa),\apa) < \tb(\Dc(\ppath),\apa)$. By Proposition \ref{rem:image_of_a}, $\dist(\xai(\apb),\apa) \le \dist(\xai(\apa),\apa)$, so it suffices to show that the defender's capture set cannot reach $\xai(\apa)$ before the attacker reaches $\apa$. 

If the attacker moves towards $\apa$ along $\spath(\xan,\apa)$ with maximum speed, then $\xai(\apa)$ will move towards $\apa$ along $\spath(\xai(\apa),\apa)$ at the same speed. Since $\ta(\xa,\apa)=\ta(\xai(\apa),\apa)<\tb(\Dc(\ppath),\apa)$, $\xa$ will reach $\apa$ before the defender capture set $\Dc(\xb,t)$ does. 

Next we show, by contradiction, that if $\xa\notin\rpa \cup \rpb$, then the attacker cannot reach $\sac$ without being captured. Suppose $\pa$ will reach some point $\pprime$ before $\Dc(\pb,t)$ does, i.e. $\dist(\xan,\pprime)<\dist(\Dc(\xbn),\pprime)=\dist(\Dc(\ppath),\pprime)$. WLOG, assume $\pprime\in\spath(\ppath,\apb)$, and note that $\dist(\Dc(\ppath),\apb)<\dist(\xan,\apb)$ since the attacker is not in $\rpb$. Starting with the definition of the shortest path, we have
\bq
\begin{aligned}
\dist(\xan,\apb) &\le \dist(\xan,\pprime) + \dist(\pprime,\apb) \\
&< \dist(\Dc(\ppath),\pprime) + \dist(\pprime,\apb) \\
&= \dist(\Dc(\ppath),\apb) \\
\dist(\xan,\apb)&< \dist(\xan, \apb) \quad \text{(since $\xan\notin\rpa$)}
\end{aligned}
\eq
This is a contradiction. Therefore, the attacker cannot cross any point $\pprime$ on $\pathd$ without being captured.
\end{IEEEproof}

If $\vela<\velb$, $\pa$ being outside of $\rpa \cup \rpb$ becomes a sufficient condition for the strong defendability of $\pathd$.

In general, $\xbn$ may not be on $\pathd$. In this case, if the defender can arrive at $\ppath$ before the attacker moves into $\rpa(\ppath)\cup\rpb(\ppath)$, then $\pathd$ is strongly defendable. Thus, given $\xjn = (\xan,\xbn)$, we may naively check whether a path $\pathd$ is strongly defendable by checking whether there exists some $\ppath\in\pathd$ such that $\tb(\xbn, \ppath)\le\ta\left(\xan,\rpa(\ppath)\cup\rpb(\ppath)\right)$. If so, then $\pathd$ is strongly defendable. The next lemma shows that it is necessary and sufficient to check whether \textit{one} special point, $\pstar\in\pathd$, can be the first arrival point for strongly defending $\pathd$.

\begin{rem} \label{rem:time_to_region_a}
Given $\ppath\in\pathd$, $\dist\left(\xan,\rpa(\ppath)\right) = \dist(\xan,\apa) - \dist(\Dc(\ppath),\apa)$. Similarly, $\dist\left(\xan,\rpb(\ppath)\right) = \dist(\xan,\apb) - \dist(\Dc(\ppath),\apb)$.
\end{rem}

\begin{lem} \label{lem:pstar}
Define $\pstar\in\pathd$ such that $\ta(\xan,\rpa)=\ta(\xan,\rpb)$. Then, $\pathd$ is strongly defendable if and only if the defender can defend $\pathd$ by first going to $\pstar$.
\end{lem}

\begin{IEEEproof}
One direction is clear by definition.

We now show the other direction's contrapositive: if the defender cannot defend $\pathd$ by first going to $\pstar$, then $\pathd$ is not strongly defendable. Equivalently, we show that if choosing $\pstar$ as the first entry point does not allow the defender to defend $\pathd$, then no other entry point does.

Suppose that the defender cannot defend $\pathd$ by choosing $\pstar$ as the first entry point, but can defend $\pathd$ by choosing another entry point $\pprime$. WLOG, assume $\dist(\Dc(\pstar),\apa)>\dist(\Dc(\pprime),\apa)$. This assumption moves $\pprime$ further away from $\apa$ than $\pstar$, causing $\rpa$ to move closer to $\xan$. Starting with Remark \ref{rem:time_to_region_a}, we have
\bq
\begin{aligned}
\dist\left(\xan,\rpa(\pstar)\right) & = \dist(\xan,\apa) - \dist(\Dc(\pstar),\apa) \\
\ta\left(\xan,\rpa(\pstar)\right) & = \ta(\xan,\apa) - \ta(\Dc(\pstar),\apa) 
\end{aligned}
\eq

Similarly, for the point $\pprime$, we have
\bq
\begin{aligned}
\dist\left(\xan,\rpa(\pprime)\right) & = \dist(\xan,\apa) - \dist(\Dc(\pprime),\apa) \\
\ta\left(\xan,\rpa(\pprime)\right) & = \ta(\xan,\apa) - \ta(\Dc(\pprime),\apa) 
\end{aligned}
\eq

Then, subtracting the above two equations, we see that the attacker can get to $\rpa$ sooner by the following amount:
\bq
\begin{aligned}
& \ta\left(\xan,\rpa(\pstar)\right) - \ta\left(\xan,\rpa(\pprime)\right) \\
&= \ta(\Dc(\pprime),\apa) - \ta(\Dc(\pstar),\apa) \\
&=\ta(\pprime,\pstar) \ge \tb(\pprime,\pstar)
\end{aligned}
\eq

We now show that the defender can get to $\pprime$ sooner than to $\pstar$ by less than the amount $\tb(\pprime,\pstar)$, and in effect ``gains less time" than the attacker does by going to $\pprime$. We assume that $\pprime$ is closer to the defender than $\pstar$ is (otherwise the defender ``loses time" by going to $\pprime$). By the triangle inequality,
\bq
\begin{aligned}
\dist(\xbn,\pstar) & \leq \dist(\xbn,\pprime) + \dist(\pprime,\pstar) \\
\dist(\xbn,\pstar) - \dist(\xbn,\pprime) & \leq  \dist(\pprime,\pstar) \\
\tb(\xbn,\pstar) - \tb(\xbn,\pprime) & \leq  \tb(\pprime,\pstar)
\end{aligned}
\eq
\end{IEEEproof}
Lemmas \ref{lem:d_winning_region} and \ref{lem:pstar} give a simple algorithm to compute, given $\xan$, the region that the defender must be in for a path of defense $\pathd$ to be strongly defendable:
\begin{enumerate}
\item Given $\apa, \apb,\xan$, compute $\pstar$ and $\rpa(\pstar),\rpb(\pstar)$.
\item If $\vela=\velb$, then $\pathd$ is strongly defendable if and only if $\xbn\in\dr(\apa,\apb;\xan)=\{x:\tb(\x,\pstar) \le \ta(\xan,\rpa \cup \rpb)\}$. 
\end{enumerate}

\begin{figure}
\centering
\includegraphics[width=0.45\textwidth]{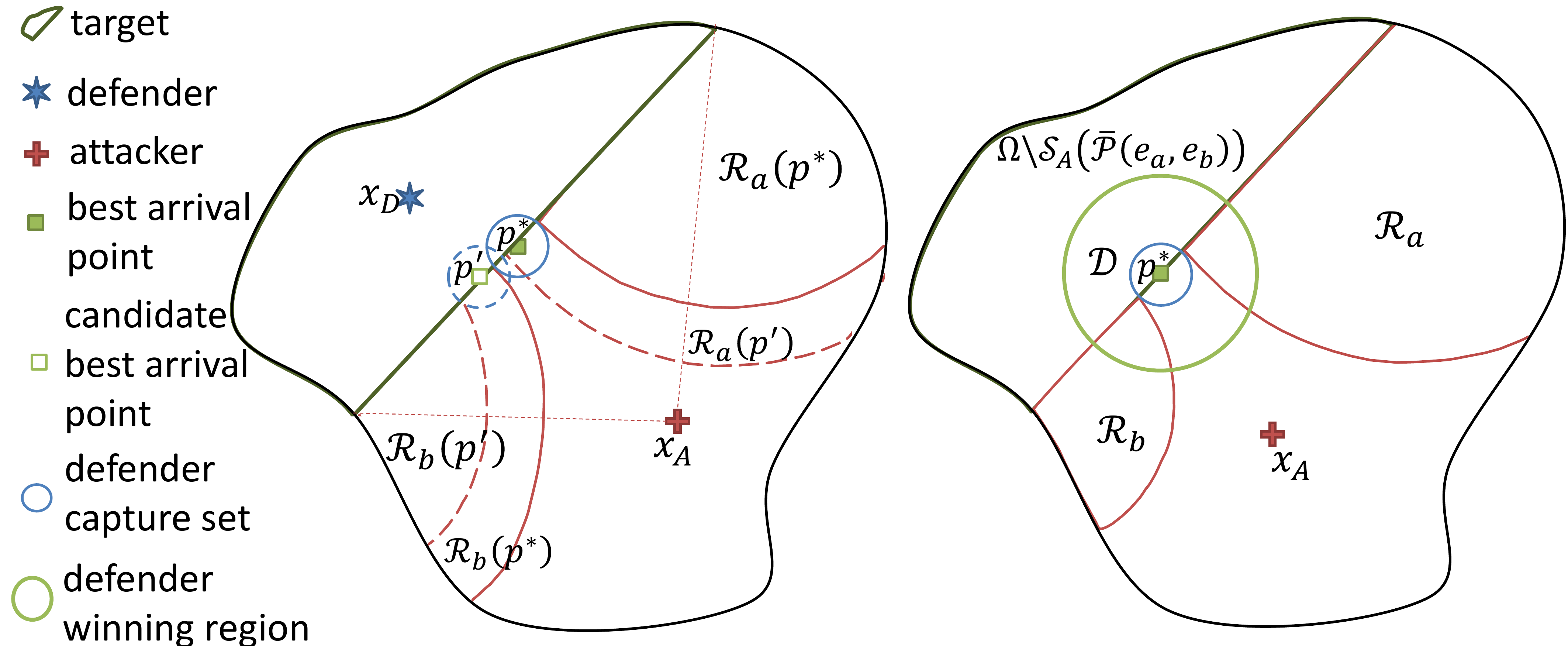}
\caption{Left: If the defender cannot defend $\pathd$ by first going to $\pstar$, then he cannot defend $\pathd$ by going to any other point $\ppath$. Right: Defender winning region $\dr$.}
\label{fig:lemma2}
\end{figure}

The computations in this algorithm can be efficiently done by solving a series of 2D Eikonal equations by using FMM \cite{Sethian1996}, reducing our 4D problem to 2D. Fig. \ref{fig:lemma2} illustrates the proof of Lemma \ref{lem:pstar} and the defender winning region $\dr$.

\subsection{The Path Defense Solution to the Reach-Avoid Game}
\label{subsec:reach_avoid}
The central idea in using path defense is that if the target set is enclosed by some strongly defendable path for some $\apa,\apb$, then the defender can win the game using the semi-open-loop strategy outlined in this section, \textit{even if} the attacker uses the optimal control. Checking for strongly defendable paths adds more conservatism towards the defender, but makes computation much more efficient.

Naively, one could fix $\apa$, then search all other anchor points $\apb\in\boundary$ to find a defendable path. However, we can reduce the number of paths that needs to be checked by only checking paths of defense $\pathd$ that touch the target set. In a simply connected domain, this reduction in the number of paths checked does not introduce any additional conservatism. 



If some strongly defendable path $\pathd$ encloses the target set, then the defender's strategy would be to first go to $\pstar\in\pathd$ (an open-loop strategy), then move towards $\xai(\apa)$ or $\xai(\apb)$ until the level set image is captured (a closed-loop strategy). Finally, the defender can simply track the captured level set image (a closed-loop strategy). This is a ``semi-open-loop" strategy. The following algorithm approximates a 2D slice conservatively towards the defender:

\begin{alg}~ \label{alg:PD_RA}
Given attacker position,
\begin{enumerate}
\item Choose some point $\apa\in\boundary$, which defines $\apb$ to create a path of defense $\pathd$ that touches the target $\target$. \label{step:createPath}
\item Repeat step \ref{step:createPath} for a desired set of points $\apa\in\boundary$. \label{step:repeatCreatePath}
\item For some particular $\pathd$, determine the defender winning region $\dr(\apa,\apb;\xan)$.\label{step:dWinRegion}
\item Repeat step \ref{step:dWinRegion} for all the paths created in steps \ref{step:createPath} and \ref{step:repeatCreatePath}.
\item The union $\bigcup_{\apa} \dr(\apa,\apb;\xan)$ gives the approximate 2D slice, representing the conservative winning region for the defender in the two-player reach-avoid game. \label{step:union}
\end{enumerate}
\end{alg}


\section{From Two-Player to Multiplayer} \label{sec:two_to_multi}

\subsection{Maximum Matching}
\label{subsec:max_match}
We piece together the outcomes of all attacker-defender pairs using maximum matching as follows:

\begin{alg}~
\begin{enumerate}
\item Construct a bipartite graph with two sets of nodes $\pas,\pbs$. Each node represents a player.
\item For all $i,j$, draw an edge between $\pbm{i}$ and $\pam{j}$ if $\pbm{i}$ wins against $\pam{j}$ in a two-player reach-avoid game.
\item Run any matching algorithm (eg. \cite{Schrjiver2004, Karpinski1998}) to find a maximum matching in the graph. 
\end{enumerate}
\end{alg}

\begin{figure}
\centering
\includegraphics[width=0.4\textwidth]{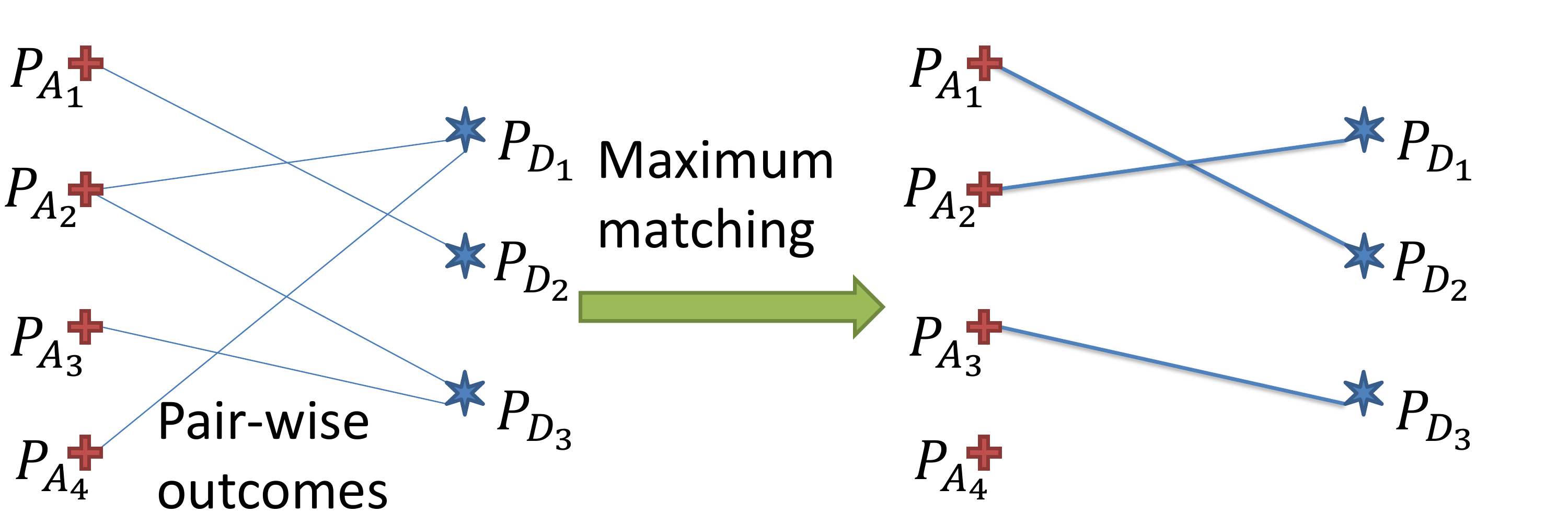}
\caption{An illustration of using maximum matching to conservatively approximate the multiplayer reach-avoid game.}
\label{fig:general_procedure}
\end{figure}

After finding a maximum matching, we can guarantee an upper bound on the number of attackers that is be able to reach the target. If the maximum matching is of size $\mm$, then the defending team would be able to prevent \textit{at least} $\mm$ attackers from reaching the target, and thus $\NA-\mm$ is an upper bound on the number of attackers that can reach the target. The maximum matching approach is illustrated in Fig. \ref{fig:general_procedure}.

\subsection{Time-Varying Defender-Attacker Pairings}
\label{subsec:tvarp}
With the next algorithm, the bipartite graph and its corresponding maximum matching can be updated, potentially in real time, as the players change positions during the game:

\begin{alg}~
\begin{enumerate}
\item Given each $\xbm{i}$ and each $\xam{j}$, determine whether $\pam{j}$ can win against $\pbm{i}$ for all $i,j$. \label{step:pairwise}
\item Assign a defender to each attacker that is part of a maximum matching.
\item For a short duration $\Delta$, apply a winning control input and compute the resulting trajectory for each defender that is part of the maximum matching. For the rest of the defenders and for all attackers, compute the trajectories assuming some (any) control function. \label{step:traj}
\item Update the player positions after the duration $\Delta$ and repeat steps \ref{step:pairwise} to \ref{step:traj}  with the new player positions.
\end{enumerate}
\end{alg}

As $\Delta\rightarrow 0$, the above procedure continuously computes a bipartite graph and its maximum matching. As long as each defender uses a winning control input against the paired-up attacker, the size of maximum matching can never decrease. 


\subsection{Application to the Two-Player HJI Solution}
\label{subusec:MMHJI}
In general, solving $\NA\ND$ 4D HJI PDEs gives us the pairwise outcomes between every attacker-defender pair. The computation time required is thus $C \NA\ND$, where $C$ is the time required to solve a single 4D HJI PDE. The pairwise outcomes can then be merged together to approximate the $\NA$ vs. $\ND$ game. In the case where each team has a single maximum speed, solving \textit{one} 4D HJI PDE would characterize all pairwise outcomes.

Since the solution to the 4D HJI PDE characterizes pairwise outcomes based on any attacker-defender joint-state, it allows for real-time updates of the maximum matching. As players move to new positions, the pairwise outcome can be updated by simply checking whether $(\xam{i}, \xbm{j})$ is in $\mathcal{RA}_\infty(R,A)$.

\subsection{Application to the Two-Player Path Defense Solution}
\label{subsec:MMPD}
To use the pairwise outcomes determined by the path defense approach for approximating the solution to the multiplayer game, we add the following step to Algorithm \ref{alg:PD_RA}: 
\begin{enumerate}
\setcounter{enumi}{5}
\item Repeat steps \ref{step:dWinRegion} to \ref{step:union} for every attacker position.
\end{enumerate}

For a given domain, set of obstacles, and target set, steps \ref{step:createPath} and \ref{step:repeatCreatePath} in Algorithm \ref{alg:PD_RA} only need to be performed once, regardless of the number of players. In step \ref{step:dWinRegion}, the speeds of defenders come in only through a single distance calculation from $\pstar$, which only needs to be done once per attacker position. Therefore, the total computation time required is on the order of $C_1 + C_2 \NA$, where $C_1$ is the time required for steps \ref{step:createPath} and \ref{step:repeatCreatePath}, $C_2$ is the time required for steps \ref{step:dWinRegion} to \ref{step:union}. 

\subsection{Defender Cooperation}
One of the strengths of the maximum matching approach is its simplicity in the way cooperation among the defenders is incorporated from pairwise outcomes. More specifically, cooperation is incorporated using the knowledge of the strategy of each teammate, and the knowledge of which attackers each teammate can win against in a 1 vs. 1 setting. 

The knowledge of the strategy of each teammate is incorporated in the following way: When the pairwise outcomes for each defender is computed, a particular defender strategy used. The strategy of each defender is then used to compute pairwise outcomes, which are used in the maximum matching process. Each defender may use the optimal closed-loop strategy given by the two-player HJI solution, the semi-open-loop strategy given by the two-player path defense solution, or even another strategy that is not described in this paper. In fact, different defenders may use a different strategy.

As already mentioned, all of the information about the strategy of each defender is used to compute the pairwise outcomes. Since each pairwise outcome specifies a winning region for the corresponding defender, each defender can be guaranteed to win against a set of attackers in a one vs. one setting. The set of attackers against which each defender can win is then used to construct the bipartite graph on which maximum matching is performed. While executing the joint defense strategy as a team, each defender simply needs to execute its \textit{pairwise} defense strategy against the attacker to which the defender is assigned. 

The maximum matching process optimally combines the information about teammates' strategies and competence to derive a joint strategy to prevent as many attackers from reaching the target as possible. The size of the maximum matching then guarantees an upper bound on the number of attackers that can reach the target. To our knowledge, no other method can synthesize a joint defender control strategy that can provide such a guarantee in a multiplayer game. 



\section{Numerical Results}
\label{sec:simulation}
We use a 4 vs. 4 example to illustrate our methods. The game is played on a square domain with obstacles. Defenders have a capture radius of $0.1$ units, and all players have the same maximum speed. Computations were done on a laptop with a Core i7-2640M processor with 4 gigabytes of memory.

\subsection{HJI Formulation}
\begin{figure}
	\centering
	\includegraphics[width=0.45\textwidth]{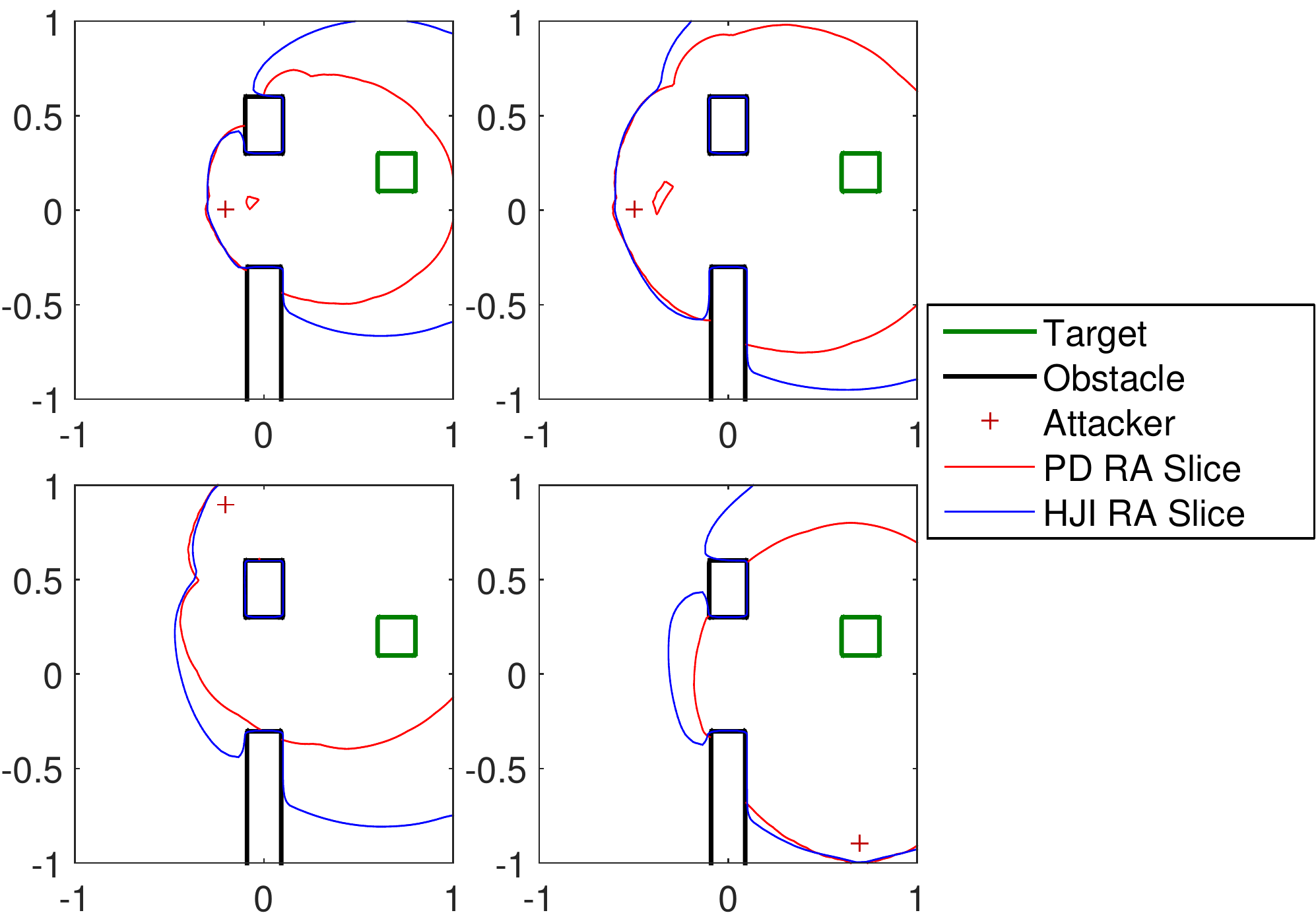}
	\caption{Reach-avoid slices computed using the HJI approach and the path defense approach. }
	\label{fig:comp_ol}
\end{figure}

Fig.  \ref{fig:comp_ol} shows the results of solving the corresponding 4D HJI PDE in blue. Computing the 4D reach-avoid set on a grid with 45 points in each dimension took approximately 30 minutes. All players have the same maximum speed, so only a single 4D HJI PDE needed to be solved. To visualize the 4D reach-avoid set, we take 2D slices of the 4D reach-avoid set sliced at each attacker's position. 


Fig. \ref{fig:mm} shows the bipartite graph and maximum matching obtained from the pairwise outcomes. The maximum matching is of size 4. This guarantees that if each defender plays optimally against the attacker matched by the maximum matching, then \textit{no} attacker will be able to reach the target.

%


\subsection{Path Defense Formulation}

Fig. \ref{fig:comp_ol} shows, in red, the results of using the path defense approach to compute conservative approximations of the 4D reach-avoid set sliced at various attacker positions.

\begin{figure}
	\centering
	\includegraphics[width=0.5\textwidth]{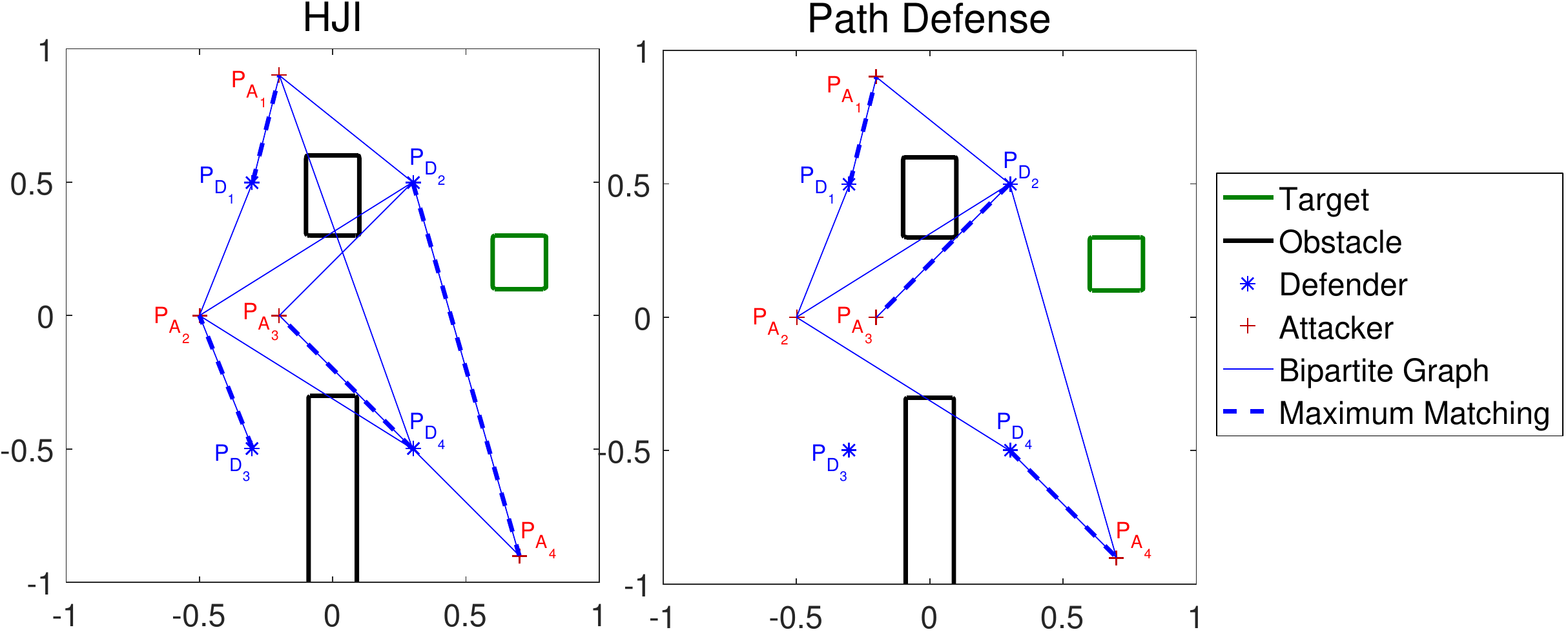}
	\caption{Merging pairwise outcomes with maximum matching.}
	\label{fig:mm}
\end{figure}

Fig. \ref{fig:mm} shows the bipartite graph and maximum matching resulting from the pairwise outcomes. In this case, the maximum matching is of size 3. This guarantees that if each defender plays against the attacker matched by the maximum matching using the semi-open-loop strategy, then \textit{at most} 1 attacker will be able to reach the target.

\begin{figure}
	\centering
	\includegraphics[width=0.35\textwidth]{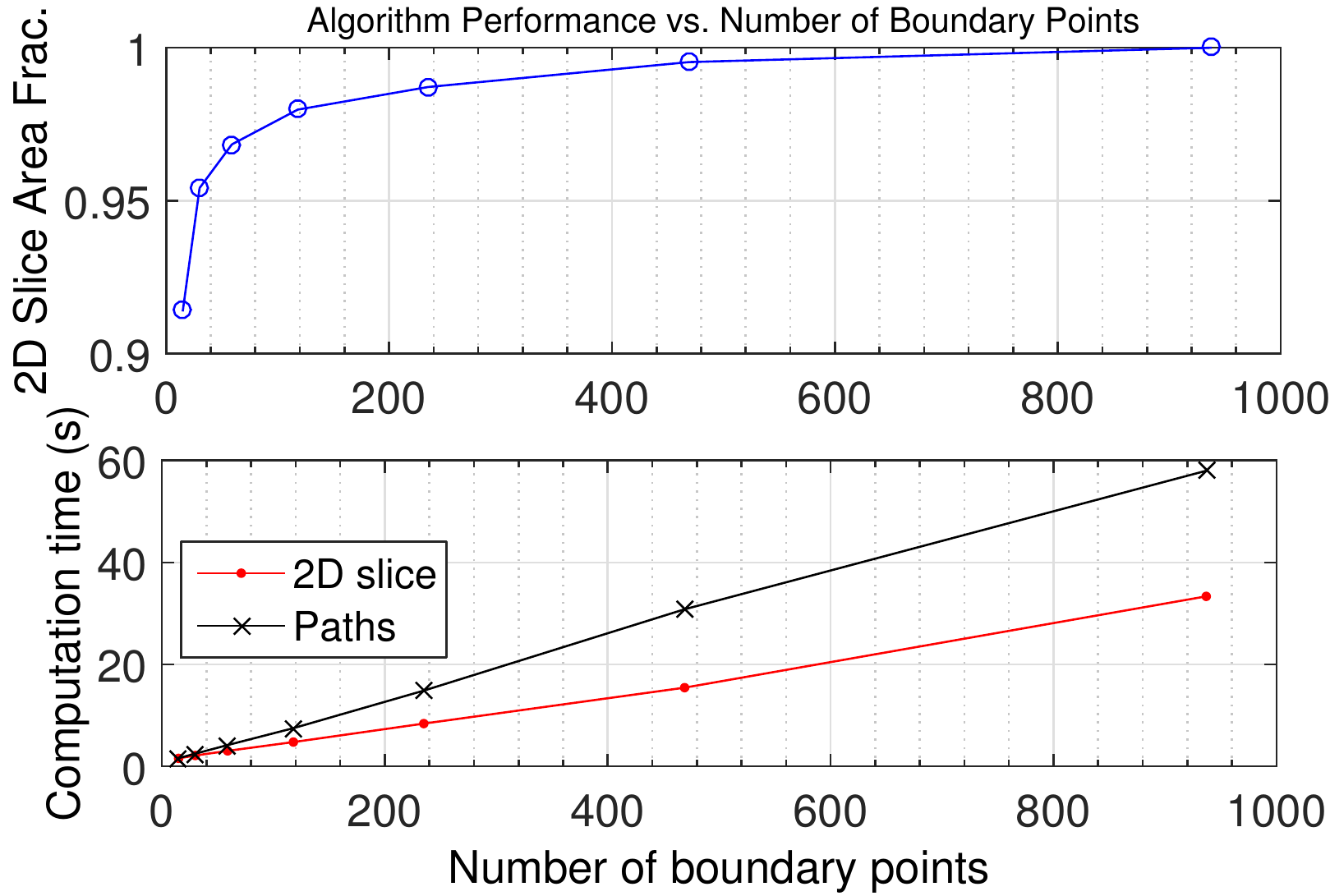}
	\caption{Performance of the Path Defense Solution.}
	\label{fig:pd_perf}
\end{figure}

Computations were done on a $200\times200$ grid, and 937 paths were used to compute the results in Fig. \ref{fig:comp_ol}. Computation time varies with the number of paths we chose in steps \ref{step:createPath} and \ref{step:repeatCreatePath} in Algorithm \ref{alg:PD_RA}. Taking the union of the defender winning regions from more paths will give a less conservative result, but requires more computation time. A summary of the performance is shown in Fig. \ref{fig:pd_perf}. With 937 paths, the computation of paths took approximately 60 seconds, and the computation of the 2D slice given the set of paths took approximately 30 seconds. However, very few paths are needed to approximate a 2D slice: Even with as few as 30 paths, the computed 2D slice covers more than 95\% of the area of the 2D slice computed using 937 paths. This reduces the computation time of the paths to 2.5 seconds, and the computation time of the 2D slices given the paths to 2.1 seconds.



\subsection{Defender Cooperation}
To highlight the element of cooperation among the defenders, we examine Fig. \ref{fig:mm} more closely. For the maximum matching results from the two-player HJI solutions (left sub-figure), we see that the maximum matching creates the pairs $(\pbm{1}, \pam{1})$, $(\pbm{2}, \pam{4})$, $(\pbm{3}, \pam{2})$, and $(\pbm{4}, \pam{3})$. In general, such a pairing is not intuitively obvious. For example, it is not obvious, without knowledge of the pairwise outcomes, that $\pbm{1}$ can win against $\pam{1}$ in a one vs. one setting. If one were not sure whether $\pbm{1}$ can win against $\pam{1}$, then $(\pbm{2}, \pam{1})$ may seem like a reasonable pair. However, if $\pbm{2}$ defends against $\pam{1}$, then$\pbm{1}$ would defend against $\pam{2}$, leaving $\pbm{3}$ unable to find an attacker that $\pbm{3}$ can be guaranteed to win against to pair up with. The same observations can be made in many of the other pairings in both sub-figures.

For the maximum matching results from the two-player path defense solutions (right sub-figure), a semi-open-loop strategy is used. In this case, given the same initial conditions of the two teams, each defender may only be guaranteed to successfully defend against fewer attackers in a one vs. one setting compared to when using the optimal closed-loop strategy from the two-player HJI solution. Regardless of the strategy that each defender uses, the maximum matching process can be applied to obtain optimal defender-attacker pairs given the strategy used and the set of attackers each defender can be guaranteed to win against in a 1 vs. 1 setting.

\subsection{Real-Time Maximum Matching Updates}
\begin{figure}
	\centering
	\includegraphics[width=0.4\textwidth]{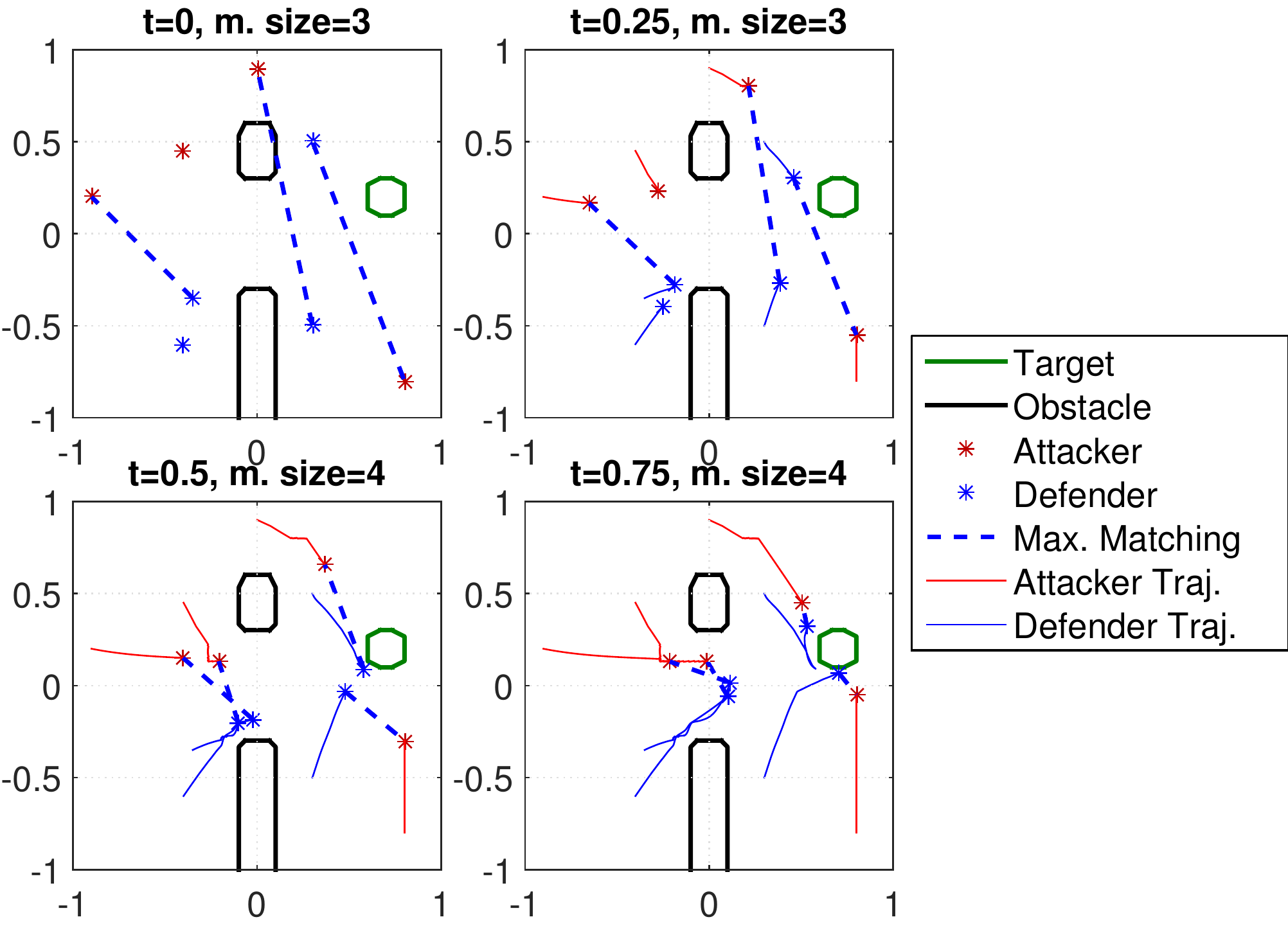}
	\caption{Increasing maximum matching size over time.}
	\label{fig:real_time_update}
\end{figure}

After determining all pairwise outcomes, pairwise outcomes of \textit{any} joint state of the attacker-defender pair are characterized by the HJI approach. This allows for updates of the bipartite graph and its maximum matching as the players play out the game in real time. Fig. \ref{fig:real_time_update} shows the maximum matching at several time snapshots of a 4 vs. 4 game. Each defender that is part of a maximum matching plays optimally against the paired-up attacker, and the remaining defender plays optimally against the closest attacker. The attackers' strategy is to move towards the target along the shortest path while steering clear of the obstacles by $0.125$ units. The maximum matching is updated every $\Delta=0.005$ seconds. At $t=0$ and $t=0.2$, the maximum matching is of size 3, which guarantees that at most one attacker will be able to reach the target. After $t=0.4$, the maximum matching size increases to 4, which guarantees that no attacker will be able to reach the target.

\section{Conclusion and Future Work}
\label{sec:conclusion}
A multiplayer reach-avoid game is numerically intractable to analyze by directly solving the corresponding high dimensional HJI PDE. To address this, we presented a way to tie together pairwise outcomes using maximum matching to approximate the solution to the full multiplayer game, guaranteeing an upper bound on the number of attackers that can reach the target. We also presented two approaches for determining the pairwise outcomes. The HJI approach is computationally more expensive, produces the optimal closed-loop control strategy for each attacker-defender pair, and efficiently allows for real time maximum updates. The path defense approach is conservative towards the defender, performs computation on the state space of a single player as opposed to the joint state space, and scales only linearly in the number of attackers.

\bibliographystyle{IEEEtran}
\bibliography{references}

\end{document}